\documentclass[10pt,a4paper]{article}

\usepackage[all]{xy}

\usepackage{epsfig}
\usepackage{amssymb}
\usepackage{amsmath}

\usepackage{graphics,graphicx}
\usepackage{eepic}

\usepackage{epsfig}
\usepackage{trees}
\usepackage{vmargin}
\usepackage{theorem}

\newtheorem{theo}{Theorem}[section]
\newtheorem{theorem}[theo]{Theorem}
\newtheorem{lemma}[theo]{Lemma}
\newtheorem{proposition}[theo]{Proposition}
\newtheorem{remark}[theo]{Remark}
\newtheorem{defin}[theo]{Definition}

{\theorembodyfont{\rmfamily}
\newtheorem{example}[theo]{Example}}

\newenvironment{proof}{\noindent\text{\textbf{Proof.\:}}}{}

\def\qed{\hfill $\square$ \goodbreak \bigskip}

\newcommand{\R}{\mathbb{R}}

\title{Structure of shape derivatives\\ around irregular domains and applications}
\author{Jimmy Lamboley and Michel Pierre\footnote{Antenne de Bretagne de l'ENS Cachan et IRMAR, Campus de Ker Lann, 35170-BRUZ, France, jimmy.lamboley@bretagne.ens-cachan.fr, michel.pierre@bretagne.ens-cachan.fr}}

\begin{document}
\maketitle
\begin{abstract}
In this paper, we describe the structure of shape derivatives
around sets which are only assumed to be of finite perimeter in
$\R^N$. This structure allows us to define a useful notion of
positivity of the shape derivative and we show it implies its
continuity with respect to the uniform norm when the boundary is
Lipschitz (this restriction is essentially optimal).  We apply
this idea to various cases including the perimeter-type
functionals for convex and pseudo-convex shapes or the Dirichlet
energy of an open set.

{\it Keywords:\,} Shape optimization, shape derivatives, sets of
finite perimeter, convex sets, Dirichlet energy
\smallskip

\end{abstract}

\section{Introduction and main results}\label{sect:intro}
The goal of this paper is to describe and to exploit {\em the structure of shape derivatives around irregular shapes}. We will only assume that they are of finite perimeter and use the tools of geometric measure theory.

Let us first recall the famous {\em Hadamard's structure theorem} around {\em regular} shapes. Let $E$ be a regular open subset of $\R^N$. We introduce small perturbations of $E$ of the form $E_\theta:=(Id+\theta)(E)$ where $Id=identity$ and $\theta$ is a small element of the space
$$\Theta:=\{\theta:\R^N\rightarrow \R^N;\; \theta \textrm{ is } C^1,\;\theta \textrm{ and }D\theta
\textrm{ are bounded}\},$$ endowed with the norm
$\|\theta\|_{1,\infty}:=\sup_{x\in\R^N}|\theta(x)|+\sup_{x\in\R^N}\|D\theta(x)\|$,
where $\|\cdot\|$ is a chosen norm for linear mappings from $\R^N$
into itself. We consider a shape functional
$$
\begin{array}{rcl} \mathcal{J}_E :  \Theta &\rightarrow &\R\\\theta&
\mapsto&\mathcal{J}_E(\theta):=J((Id+\theta)(E)),
\end{array}
$$
defined at least in a neighborhood of the origin, where $J:\{E_\theta\subset \R^N;\,\theta \;small\}\to \R$ is given. We denote by $\mathcal{J}_E'(0)$ the Fr\'echet derivative (or differential) of  $\mathcal{J}_{E}$ at $\theta=0$
 when it exists. Then, Hadamard's structure theorem states that this derivative depends only on the normal component of the deformations at the boundary of $E$. More precisely (see the pioneer work \cite{Had}, or \cite{MS},
\cite{DZ} p. 348, \cite{HP} p. 219, and \cite{NP}):
\begin{proposition}\label{prop1}
Let $E$ be an open subset of $\R^N$ of class $\mathcal{C}^2$ and assume that $\mathcal{J}_E$ is (Fr\'echet) differentiable at $\theta=0$ in $\Theta$. Denote by $\nu_E$ the outward unit normal
derivative to $\partial E$. Then, there exists $l
:\mathcal{C}^1(\partial E,\R)\rightarrow \R$ a continuous linear form
such that
\begin{eqnarray}\label{struct1}
\forall \xi \in \Theta, \;\;\mathcal{J}_E'(0)\cdot\xi =
l(\xi_{{|\partial E}}\cdot\nu_E),
\end{eqnarray}
where $\xi_{|\partial E}$ denotes the restriction of $\xi$ to the
topological boundary $\partial E$ of $E$.
\end{proposition}

Without any assumption on the regularity of $E$, we at least have the following (see \cite{DZ} p. 348 or \cite{HP} p. 217):
\begin{proposition}\label{prop2}
Let $E\subset \R^N$. We assume that $\mathcal{J}_E$
is differentiable at 0 in $\Theta$.
Then $\mathcal{J}_E'(0)$ is a distribution in $\R^N$ of order at most 1 and whose support is included in $\partial E$.\\
\end{proposition}

A first question we address here is the following: {\em what is left of the structure (\ref{struct1}) when $E$ is irregular?}\\

Throughout the paper, we will assume that the functional $\mathcal{J}_E$ satisfies
\begin{eqnarray}\label{hypJ}
{\rm for \;all\;small}\;\theta_1,\theta_2\in\Theta,\;\;[E_{\theta_1}=E_{\theta_2}\;a.e.]\;\;\Rightarrow \;\;[\mathcal{J}_E(\theta_1)=\mathcal{J}_E(\theta_2)]\,.
\end{eqnarray}
As $Id+\theta$ is a diffeomorphism for $\|\theta\|_{1,\infty}$ small, this is obviously satisfied as soon as $E$ is regular enough since then: $E_{\theta_1}=E_{\theta_2}\;a.e.$ implies $E_{\theta_1}=E_{\theta_2}$ "everywhere". When $E$ is irregular, it is also satisfied for many functionals: it actually means that $J$ (or $\mathcal{J}_E$) is a function of the characteristic function of $E_\theta$, seen as a class of functions defined a.e.. This contains even $H^1$-energy functionals (see examples in Section 3).

We will also assume that {\em $E$ is of finite perimeter} so that we may use the tools of geometric measure theory (see the Appendix). In particular, we know that the unit normal $\nu_E$ may be defined everywhere on the so-called {\em reduced boundary} $\Gamma:=\partial^*\! E$. Moreover, if $\mathcal{H}^{N-1}$ denotes the $(N-1)$-dimensional Hausdorff measure (see the Appendix), the restriction mapping
$$\Phi: \xi\in\Theta\mapsto \xi_{|\Gamma}\cdot \nu_E\in L^\infty(\Gamma,\mathcal{H}^{N-1})$$
is well-defined. We then prove the following first main result (where $\Phi(\Theta)$ denotes the image of $\Phi$):

\begin{theorem}\label{the1} Let $E$ be a set of finite perimeter and
$\Gamma=\partial^*\!E$ its reduced boundary. Then, there exists a
Banach space structure on $\Phi(\Theta)$ such that, for any
functional $\mathcal{J}_E$ satisfying (\ref{hypJ}) and
differentiable at $\theta=0$, we have
\begin{eqnarray}\label{struct2}
\forall \xi\, \in \Theta, \;\;\mathcal{J}_E'(0)\cdot\xi =
l(\xi_{|\Gamma}\cdot \nu_E),
\end{eqnarray}
where $l$ is a continuous linear form on $\Phi(\Theta)$.
\end{theorem}

When $E$ is regular as in Proposition \ref{prop1}, $\Phi(\Theta)$ is nothing but $\mathcal C^1(\partial\Omega)$ endowed with its $C^1$-norm. In all cases, it is stronger than or equal to the uniform norm ($=L^\infty(\Gamma)$-norm). In general, $\Phi(\Theta)$ and its norm are difficult to identify explicitly.  We try to do it next in some specific cases. In particular, we are interested in cases where the linear form $l$ is actually {\em continuous with respect to the uniform norm}, that is when $l$ may essentially be considered as a Radon measure on $\overline{\Gamma}$.

In this context, a next natural question is then: {\em Does positivity imply continuity of the linear form $l(\cdot )$ with respect to the uniform norm?}

Thanks to the structure Theorem \ref{the1}, we may indeed define a natural notion of positivity for  $\mathcal{J}_E'(0)$ by saying that $l$ is positive: for all $\xi\in\Theta$
\begin{eqnarray}\label{pos}
[\;\xi_{|\Gamma}\cdot \nu_E\geq 0\;\;\mathcal{H}^{N-1}-a.e.\;on\;\Gamma] \;\Rightarrow \;[\;\mathcal{J}_E'(0)\cdot\xi=l(\xi_{|\Gamma}\cdot \nu_E)\geq 0]\;.
\end{eqnarray}

As we will see, this property will simply be satisfied as soon as $J$ is nondecreasing with respect to the inclusion a.e..

A second main result of this paper is the following (see Propositions \ref{prop:pos} and \ref{Liprop}): {\em Positivity of $l$ does imply its $L^\infty$-continuity when $E$ is an open set with a Lipschitz boundary}. And this is rather optimal: we provide an example showing that it is not the case as soon as $\partial E$ has a singularity like a cusp. Surprisingly enough, it is then possible that, at the same time, $\mathcal{J}_E'(0)$ is $L^\infty$-continuous and $l$ is not. \\

The rest of the paper is devoted to applications of the
$L^\infty$-continuity result. For instance, it is known that, if
$\mathcal{P}_E(\theta):=P((Id+\theta)(E))$ where $P(\cdot )$ is
the perimeter function and $P(E)<+\infty$, then
$\mathcal{P}_E(\cdot )$ is differentiable at $\theta=0$ (for {\bf
any} set with $P(E)<+\infty$) and its derivative is given by
\begin{eqnarray}\label{per}
\forall \xi\in\Theta,\;\;\mathcal{P}'_E(0)\cdot\xi=\int_\Gamma \nabla_\Gamma\cdot \xi\,d\mathcal{H}^{N-1},
\end{eqnarray}
where $\nabla_\Gamma\cdot \xi:=\nabla\cdot
\xi-(D\xi\cdot\nu_E)\cdot \nu_E$ denotes the tangential divergence
of $\xi$ on $\Gamma$. If $E$ is regular enough, we may rewrite
this as
\begin{eqnarray}\label{per1}
\forall \xi\in\Theta,\;\;\mathcal{P}'_E(0)\cdot\xi=\int_\Gamma H_\Gamma\,\xi\cdot \nu_E\,d\mathcal{H}^{N-1},
\end{eqnarray}
where $H_\Gamma$ denotes the mean curvature of $\Gamma=\partial
E$. Thus, Theorem \ref{the1} provides a way to define -in a very
weak sense- the mean curvature of any set with finite perimeter
(see also \cite{BGM} and \cite{BGT}). Next, our approach of
positivity, leads to more precise results for convex sets. For
instance, we state:

\begin{theorem}\label{thm:convex}
Let $C$ be a bounded convex subset of $\R^N$. Then there exists $k>0$ such that
\begin{equation}\label{PC1}
\forall \xi \in \Theta,\;\; |\mathcal{P}_C'(0)\cdot\xi|\; \leq
\;k\,\|\xi_{|\partial C}\cdot\nu_C\|_{L^\infty(\partial
C,\mathcal{H}^{N-1})}\;.
\end{equation}
In particular, there exists
an $\R^N$-valued measure $\overrightarrow{\textbf{H}}$ with support in $\partial C$ such that
\begin{equation}\label{PC2}
\forall \xi \in \Theta, \;\;\mathcal{P}_C'(0)\cdot\xi =
\int_{\partial C} \nabla_{\!\partial C}\cdot  \xi \,
d\mathcal{H}^{N-1} = \int_{\partial C}\xi\cdot
d\overrightarrow{\textbf{H}}\,.
\end{equation}
\end{theorem}

\begin{remark}\rm
The existence of the vector-valued {\em mean curvature measure} $\overrightarrow{\textbf{H}}$ for a convex set is natural. A proof of its existence may be found for instance in \cite{BGM} (through a quite different way). Here, we obtain it as representing the shape derivative of the perimeter; moreover, its Radon-measure property is obtained by a rather elementary approach relying on the existence and the positivity of the linear mapping $l(\cdot )$:  note that the second part of the theorem is indeed an immediate consequence of the first one since it implies
$$\forall \xi \in \Theta,\;\; |\mathcal{P}_C'(0)\cdot\xi| \leq
k\,\|\xi\|_{L^\infty(\partial C,{\mathcal H}^{N-1})}\leq k\|\xi\|_{L^\infty(\R^N,\R^N)},$$
from which the existence of $\overrightarrow{\textbf{H}}$ follows. But, the estimate (\ref{PC1}) is quite more precise than the only existence of $\overrightarrow{\textbf{H}}$.
\end{remark}

 We also deduce from the above considerations that the derivative of $\theta \mapsto \int_{\partial C_\theta} g$ (where $C_\theta=(Id+\theta)(C)$) is also a measure for regular enough functions $g$.\\

 As stated in Proposition \ref{pseudo}, properties of Theorem \ref{thm:convex} may be extended to Lipschitz {\em pseudo-convex sets} and this improves a result in \cite{BGM}.\\

 Finally, we apply our approach to quite different classes of functionals appearing in shape optimization problems associated with partial differential equations. We do it, in particular, for the Dirichlet energy of a variable open (or even measurable) bounded subset of $\R^d$, namely
 $$\theta\mapsto \mathcal{J}_E(\theta):=\int_{E_\theta}\frac{1}{2}|\nabla u_\theta|^2-f\,u_\theta,$$
 where $f:\R^N\to \R$ is given and $u_\theta$ is the solution of the Dirichlet problem
 $$-\Delta u_\theta=f\;\textrm{ on } E_\theta=(Id+\theta)(E),\;\;u_\theta=0\;\textrm{ on }\partial E_\theta.$$
 We also consider the functionals $\theta\mapsto\lambda_k(E_\theta)$ where $\lambda_k(\cdot )$ denotes the $k$-th eigenvalue of the Laplace operator with homogeneous Dirichlet boundary conditions.\\

\section{Proof of the structure theorem and examples}
\label{sect:Structure}

\noindent{\bf Proof of Theorem \ref{the1}:}\\

\textbf{Preliminary remark.}
Let us recall the two main steps of a proof in the case $E$ is
regular. First,
$$[\xi\cdot \nu_E \equiv 0\; on \;\Gamma]\implies [\mathcal{J}_E'(0)\cdot\xi =0].\;$$
(This comes from the invariance of $\Gamma$ by the flow of $\xi$).
Then, this allows one to factorize $\mathcal{J}_E'(0)$ through the
quotient space $\Theta/K$ where $K:=\{\xi\in\Theta; \xi\cdot
\nu_E=0\;on\;\Gamma\}$. The second step of the proof is to
identify this quotient space as being isometric to
$\mathcal{C}^1(\Gamma)$ by studying the mapping
$$
\begin{array}{cccc} \Phi : &\Theta &\rightarrow &\mathcal{C}^1(\Gamma)
\\&\xi&
\mapsto&\xi_{|\Gamma}\cdot \nu_E\;.
\end{array}
$$
One can see that $\Phi$ induces a isomorphism $\widetilde\Phi$ between
$\Theta/K$ and $\mathcal{C}^{1}(\Gamma)$
for their usual norms, with the use of a continuous extension operator from $\mathcal{C}^1(\Gamma)$ to $\Theta$.

These two steps may be reproduced in the general case as follows.\\

\textbf{Step 1 :} Let $\xi$ in $\Theta$ such that $\xi\cdot
\nu_E=0$ $\mathcal{H}^{N-1}$-a.e. on $\Gamma$. Let $\gamma_t$ be
its associated flow (see the statement of Lemma \ref{deriv} in the
Appendix). Let us prove that $\gamma_t(E)=E$ a.e. for all $t$. It
will follow from assumption (\ref{hypJ}) that
$\mathcal{J}_E(\gamma_t-I)=J(\gamma_t(E))=J(E)$ for all $t$ and
consequently $\mathcal{J}_E'(0)\cdot\xi=0$ since then
$$0=\frac{d}{dt}_{|t=0}\mathcal{J}_E(\gamma_t-I)=\mathcal{J}_E'(0)\cdot \frac{d}{dt}_{|t=0}\gamma_t=\mathcal{J}'_E(0)\cdot \xi\,. $$

From Lemma \ref{deriv}, we have that $\frac{\partial}{\partial t}(\chi_E\circ\gamma_t)=0$ in the sense of distributions in $\R\times\R^N$ since, for all $\varphi\in \mathcal{C}_0^\infty(\R\times\R^N)$:
$$\int_{\R\times\R^N}\frac{\partial\varphi}{\partial t}\chi_E\circ\gamma_t=\int_\R\int_{\partial^*\!E}(\varphi\circ\gamma_t^{-1})\det D(\gamma_t^{-1})(\xi\cdot\nu_E)\,d\,\mathcal{H}^{N-1}dt=0.$$

Since $t\to\chi_E\circ\gamma_t\in L^1_{loc}(\R^N)$ is continuous, it follows that $\chi_E\circ\gamma_t=\chi_E$ a.e. for all $t$.\\

\textbf{Step 2 :} Thanks to Step 1, we can factorize $\mathcal{J}_E'(0)$ as:
$$
  \xymatrix{
    \Theta \ar[r]^{\mathcal{J}_E'(0)} \ar[d]_\pi  & \R \\
    \Theta/K \ar[ru]^{\widetilde{l}} &
  }
$$
where $K:=\{\xi \in \Theta, \xi\cdot \nu_E=0\;\;\mathcal{H}^{N-1}-a.e.\;on\;\Gamma\}$ and $\pi$ is the canonical projection on the quotient space. Next, the mapping
$$
\begin{array}{cccc} \Phi : &\Theta &\rightarrow &L^\infty(\Gamma,\mathcal{H}^{N-1})
\\&\xi&
\mapsto&\xi_{|\Gamma}\cdot \nu_E\;,
\end{array}
$$
induces the (algebraic) isomorphism $\widetilde{\Phi}$
$$
  \xymatrix{
    \Theta \ar[r]^{\Phi} \ar[d]_\pi  & \Phi(\Theta) \\
    \Theta/K \ar[ru]^{\widetilde{\Phi}} &
  }
$$
Since $K$ is a closed subspace of the Banach space $\Theta$, then $\Theta/K$ is also a Banach space for the induced quotient norm. We now choose to equip the image $\Phi(\Theta)=\widetilde{\Phi}(\Theta/K)$ with the transported norm, namely
$$\forall \varphi \in \Phi(\Theta), \;\|\varphi\|:=\|\widetilde{\Phi}^{-1}(\varphi)\|_{\Theta/K}=\inf\{\|\xi\|_{1,\infty}; \xi\cdot \nu_E=\varphi\;\; \mathcal{H}^{N-1}-a.e.\},$$
which provides a Banach space structure for the space $\Phi(\Theta)$ and $\widetilde{\Phi}$ is an isometry. This is summarized in the following diagram:
$$
  \xymatrix{
    \Theta \ar[rr]^{\mathcal{J}_E'(0)} \ar[d]^\pi \ar@/_1pc/[dd]_\Phi & & \R \\
    \Theta/K \ar[rru]^{\widetilde{l}} \ar[d]^{\widetilde{\Phi}} &&\\
    \Phi(\Theta)\ar[rruu]_l&&
  }
$$
We now introduce $l:=\widetilde{l}\circ\widetilde{\Phi}^{-1}$ which is a linear continuous form on the space $\Phi(\Theta)$: we have for all $\xi\in\Theta$
$$\mathcal{J}_E'(0)\cdot\xi = \widetilde{l}(\pi(\xi))=\widetilde{l}\circ\widetilde{\Phi}^{-1}(\widetilde{\Phi}\circ\pi(\xi))
=l(\Phi(\xi))=l(\xi_{|\Gamma}\cdot \nu_E).$$\qed

\begin{remark}\label{stronger}\rm The norm on $\Phi(\Theta)$ is stronger than (or equal to) the $L^\infty(\Gamma,\mathcal{H}^{N-1})$-norm. Indeed, since for $\xi\in\Theta$, $\|\xi\cdot \nu_E\|_\infty:=\|\xi\cdot \nu_E\|_{L^\infty(\Gamma,\mathcal{H}^{N-1})}\leq \|\xi\|_\infty\leq\|\xi\|_{1,\infty}$, we have for all $\varphi\in\Phi(\Theta)$
$$\|\varphi\|=\inf\{\|\xi\|_{1,\infty}; \xi\cdot \nu_E=\varphi\; \mathcal{H}^{N-1}-a.e.\}\geq \inf\{\|\xi\cdot \nu_E\|_\infty;\xi\cdot \nu_E=\varphi\;\mathcal{H}^{N-1}-a.e.\;\}=\|\varphi\|_\infty.$$
\end{remark}

\begin{remark}\rm It may be interesting to compare our invariance condition
\begin{eqnarray}\label{normale}
\xi\cdot\nu_E=0\;\mathcal{H}^{N-1}\; a.e.\; on \;\partial^*\!E,
\end{eqnarray}
 with usual Nagumo type conditions for the flow $\gamma_t$ associated with $\xi$ according to the appendix (see for instance \cite{DZ1}, \cite{DZ}). One knows that, for a closed set $F$,
 \begin{eqnarray}\label{nagumo}
 [\forall t, \gamma_t(F)\subset F] \Leftrightarrow  [\forall x\in \partial F,\;  \liminf_{\lambda \rightarrow 0} \lambda^{-1}d(x+\lambda\xi(x),F)=0].
 \end{eqnarray}
Here, the distance function is not well adapted since our sets are
defined only a.e. However, if $E\subset\R^N$ is a measurable set,
it is classical (see e.g. \cite{Giusti}) to introduce an adequate
representation of $E$ as $\widetilde{E}=(E\cup E_1)\setminus E_0$
where $E_0,E_1$ are open sets defined as follows where
$\mathcal{L}^N$ denotes the Lebesgue measure in $\R^N$ :
$$E_0=\{x\in\R^N;\;\exists r,\;\mathcal{L}^N(E\cap B(x,r))=0\},\;E_1=\{x\in\R^N;\;\exists r,\;\mathcal{L}^N(E\cap B(x,r))=\mathcal{L}^N(B(x,r))\}.$$
Then $\widetilde{E}=E$ a.e. and $\overline{\partial^*\!E}=\partial
\widetilde{E}$. By our analysis, condition (\ref{normale}) implies
that $\gamma_t(\widetilde{E})=\widetilde{E}$ a.e.. Since
$\gamma_t$ is a $C^1$-diffeomorphism, it is easy to check that this
implies the invariance of the three sets $E_1, \R^N\setminus E_0$
and $\partial\widetilde{E}$. According to (\ref{nagumo}), this is
equivalent to
$$\forall x\in \partial\widetilde{E},\; \liminf_{\lambda\to 0}\lambda^{-1}d(x+\lambda\xi(x),\partial\widetilde{E})=0.$$

Conversely, this condition does not seem
sufficient in general to imply the a.e. invariance of $E$ (or of $\widetilde{E}$) by $\gamma_t$ (or to imply (\ref{normale}) which is equivalent by Lemma \ref{deriv}). It is however the case if $E$ (or $\widetilde{E}$) is equal a.e. to an open set or to a closed set. Indeed, in this cases, we would have respectively $E=E_1$ a.e. or $E=\R^N\setminus E_0$ a.e.. It is in particular the case if $\mathcal{L}^N(\partial{\widetilde E})=0$.\\
\end{remark}

\begin{example}\label{ex} The structure proved in Theorem \ref{the1} may sometimes require some transformations to become more explicit. Let us for instance consider the perimeter function $P(\cdot )$ around
the square $[0,1]^2\subset \R^2$. An easy computation leads to the following for $\xi\in\Theta$ with compact support around the origin:
\begin{eqnarray}\label{P} \mathcal{P}_E'(0)\cdot\xi=
-\xi_1(0,0)-\xi_2(0,0)=\delta_{(0,0)}(\xi)\cdot
\left(\begin{array}{c} -1\\-1
\end{array}\right ).
\end{eqnarray}
The
dependence on $\xi\cdot \nu_E$ is rather surprising since $\nu_E$ is not defined at $(0,0)$ and $\mathcal{P}_E'(0)$ is supported in $(0,0)$.
Nevertheless, we can find explicitly the linear form :
$$\forall \varphi\in \Phi(\Theta),\; l(\varphi)=-\lim_{y\rightarrow 0} \varphi (0,y)-\lim_{x\rightarrow 0} \varphi (x,0). $$
Indeed, $-\xi_1(0,0)=\lim_{y\to 0} \xi\cdot  \nu_E(0,y)$ since
$\nu_E(0,y)=(-1,0)$ when $0<y<1$, and similarly for $\xi_2(0,0)$.
We see that $l$ is here continuous with respect to the
$L^\infty$-norm (and for the $\Phi(\Theta)$-norm as well by Remark
\ref{stronger}).

The same computation will work for any set with an angle at a point $x_0$ of its boundary. Indeed, the limits $(\nu^-,\nu^+)$ of $\nu_E(\cdot )$ from the left and from the right at $x_0$ form a basis of $\R^2$: then we can easily reconstruct the contribution to the form
 $l$ of any Dirac mass at $x_0$.\\
\end{example}
\begin{example}\label{ex0}
The situation is quite different when $E$ has a {\em cusp} at its boundary since then $(\nu^-,\nu^+)$ is not a basis any more. Let us consider for instance what happens when $E$ is, locally around the origin, the set above the graph of $y=f(x)=|x|^\alpha$ in $\R^2$ ($\alpha\in (0,1)$).

An easy computation leads, for any $\xi\in\Theta$ supported around the origin (see the next example)
$$ \mathcal{P}_E'(0)\cdot\xi=\int_\Gamma H_\Gamma \; \xi\cdot \nu_E \, d\mathcal{H}^1
-2\,\xi_2(0,0),$$ where $H_\Gamma$ is the mean curvature, defined
and integrable on the regular part $\partial E \setminus
\{(0,0)\}$. In this case, we may prove that
$\mathcal{P}_E'(0)\cdot\xi=l(\xi\cdot \nu_E)$ with
$$\forall\varphi\in\Phi(\Theta),\; l(\varphi)=\int_\Gamma H_\Gamma
\; \varphi \, d\mathcal{H}^1-\lim_{x\rightarrow 0^+} \left[\varphi
(x,|x|^\alpha) + \varphi (-x,|x|^\alpha)\right]\alpha
|x|^{\alpha-1}.$$ Indeed, for $x\neq 0$
$$[\xi\cdot \nu_E](x,f(x))=[\xi_1(x,f(x))f'(x)-\xi_2(x,f(x))]/\sqrt{1+(f'(x))^2},$$
so that
$$\lim_{x\to 0^+}[\xi\cdot \nu_E](x,x^\alpha)=\xi_1(0,0),\;\lim_{x\to 0^-}[\xi\cdot \nu_E](x,|x|^\alpha)=-\xi_1(0,0),$$
while, using that $\xi\in\mathcal{C}^1$
\begin{eqnarray}\label{xi2}
\lim_{x\to 0^+}|x|^{\alpha-1}\big[(\xi\cdot \nu_E)(x,|x|^\alpha)+(\xi\cdot\nu_E)(-x,|x|^\alpha)\big]=-2\alpha^{-1}\xi_2(0,0)\;.
\end{eqnarray}
\end{example}

\begin{example}\label{ex1}
{\bf Piecewise regular boundaries:} We may generalize the previous
computations (Examples \ref{ex}-\ref{ex0}) to any open set $E$
bounded by a closed "curve" which is piecewise regular in the
following sense: the positively oriented boundary $\partial E$ is
parametrized by its length parameter and is the image of a
continuous mapping $\zeta:[0,L]\to\R^2$ such that for some
subdivision $0=s_0<s_1<...<s_p=L$, and for all $i=1,...,p$
$$\zeta\in \mathcal{C}^1([s_{i-1},s_i])\cap \mathcal{W}^{2,1}(s_{i-1},s_i),\;H_\Gamma\in L^1(\zeta(s_{i-1},s_i)).$$
Moreover, $\zeta$ is injective on $]0,L]$ and $\zeta(0)=\zeta(L)$.

Then,
$[\xi\in\Theta\to T(\xi):=\mathcal{P}'_E(0)\cdot\xi-\int_\Gamma H_\Gamma\,\xi\cdot \nu_E\,d\mathcal{H}^1]$
is a distribution, of order at most 1, supported in the singular points. In order to identify it, assume $\xi$ is compactly supported around $\zeta_i=\zeta(s_i)$. Then, from the expression of $P((Id+t\xi)(E)$, we easily deduce
$$t^{-1}[P((Id+t\xi)(E))-P(E)]=\epsilon(t)+\int_{s_{i-1}}^{s_{i+1}}[x'(s)\frac{d}{ds}\xi_1(\zeta(s))+y'(s)\frac{d}{ds}\xi_2(\zeta(s))]\,ds.$$
Integrating by parts on $(s_{i-1},s_i)$ and $(s_i,s_{i+1})$, and
after subtracting the "regular" part, we obtain
$$T(\xi)=[x'(s_i^-)-x'(s_i^+)]\xi_1(\zeta_i)+[y'(s_i^-)-y'(s_i^+)]\xi_2(\zeta_i)=\xi(\zeta_i)\cdot[\zeta'(s_i^-)-\zeta'(s_i^+)].$$
 This be may summarized as follows:
\begin{proposition} With the above notations and assumptions
$$\forall \xi\in\Theta,\;\mathcal{P}'_E(0)\cdot\xi=\int_\Gamma H_\Gamma(\xi\cdot\nu_E)\,d\mathcal{H}^1+\sum_i\xi(\zeta_i)\cdot[\zeta'(s_i^-)-\zeta'(s_i^+)].$$
\end{proposition}
The last term in the expression above depends on the {\em
tangential components} of $\xi$ at the singular points. Although
not obvious, by Theorem \ref{the1} it is also a function of
$\xi_{|\partial^*\!E}\cdot\nu_E$. It may be made explicit by
computations as in Examples \ref{ex}-\ref{ex0}. We refer to
\cite{Glau} for similar computations and to \cite{Laurain} for
differentiation around cracks, where tangential components also
appear.
\end{example}

\section{Using the positivity of the shape derivatives}
\label{sect:Positivity}
\subsection{Positivity implies continuity...in good enough situations}
It is classical that positivity of linear forms often implies
their continuity, in particular for $L^\infty$-norms. In this
subsection, we analyze this kind of continuity properties for
shape derivatives of functionals satisfying (\ref{hypJ}). The
notion of positivity we consider for $\mathcal{J}_E'(0)$ is
defined in (\ref{pos}). It is satisfied whenever $J$ is nondecreasing with respect to the inclusion a.e.. Indeed
\begin{lemma}\label{croiss} Let $E\subset\R^N$ with finite perimeter and let $\mathcal{J}_E$ be a functional satisfying (\ref{hypJ}) and differentiable at $\theta=0$. If for all small $\theta\in\Theta$
$$ [\;E\subset E_\theta\;a.e.\;]\;\;\Rightarrow \;\;[\;\mathcal{J}_E(0)\leq\mathcal{J}_E(\theta)\;]\;,$$
then $\mathcal{J}_E'(0)\geq 0$ in the sense of (\ref{pos}).
\end{lemma}

\noindent {\bf Proof.} Let $\xi\in\Theta$ with $\xi_{|\Gamma}\cdot \nu_E\geq 0\;\;\mathcal{H}^{N-1}-a.e.\;on\;\Gamma$ and let $\gamma_t$ be its flow. Thanks to Lemma \ref{deriv}, for all $\varphi\in\mathcal{C}_0^\infty(\R\times\R^N)$ with $\varphi\geq 0$,
$$\int_{\R\times\R^N}\frac{\partial \varphi}{\partial t}\chi_E\circ\gamma_t=\int_\R\int_{\partial^*\!E}(\varphi\circ\gamma_t^{-1})\det D(\gamma_t^{-1})(\xi\cdot\nu_E)\,d\mathcal{H}^{N-1}dt\geq 0.$$
We deduce that $\frac{\partial}{\partial t}
\chi_E\circ\gamma_t\leq 0$, it follows that, for all $t\geq 0$,
$\chi_E\circ\gamma_t\leq \chi_E$ a.e. or also $E\subset
\gamma_t(E)$ a.e.. By monotonicity of $\mathcal{J}_E$, we have
$\mathcal{J}_E(0)\leq\mathcal{J}_E(\gamma_t-Id)$. Differentiating
at $t=0$ leads to  $\mathcal{J}_E'(0)\cdot \xi\geq 0$ (see Step 1
in the proof of Theorem \ref{the1}).\qed

The main estimate coming from the positivity of ${\mathcal J}_E'(0)$ is based on the following "abstract" result.

\begin{proposition}\label{prop:pos}
Let $E\subset\R^N$ be a set of finite perimeter such that
\begin{equation}\label{H}
\exists\,\xi_0 \in \Theta\;{\rm with}\;\;
{\xi_0}_{|\partial^*\!E}\cdot\nu_E \geq
\eta\;\;\;\mathcal{H}^{N-1}-\textrm{a.e.\;on }\partial^*\!E\;{\rm
for\;some}\;\eta>0\;.
\end{equation}
If a shape derivative $\mathcal{J}_E'(0)$ is positive in the sense
of (\ref{pos}), then there exists some $k\geq 0$ such that
$$\forall\, \xi \in \Theta, \;\;|\mathcal{J}_E'(0)\cdot\xi| \leq k\,\|\xi_{|\partial^*\!E}\cdot\nu_E\|_{L^\infty(\partial^*\!E,\mathcal{H}^{N-1})}\;\;.$$
In particular, the distribution $\mathcal{J}_E'(0)$ is of order 0,
that is to say, a vector valued measure (supported on $\partial
E$).
\end{proposition}

\begin{proof}
Let $\xi\in\Theta$, $\Gamma=\partial^*\!E$, $\|\xi\cdot
\nu_E\|_\infty:=\|\xi\cdot
\nu_E\|_{L^\infty(\Gamma,\mathcal{H}^{N-1})},
\zeta_0:={\xi_0}_{{|\Gamma}}\cdot\nu_E$.  We may write
$$-\|\xi\cdot \nu_E\|_\infty\frac{\zeta_0}{\eta}\leq \xi_{|\Gamma}\cdot \nu_E\leq\|\xi\cdot \nu_E\|_\infty\frac{\zeta_0}{\eta}.$$
By positivity and linearity of $l$, we deduce
$$-\|\xi\cdot \nu_E\|_\infty\frac{l(\zeta_0)}{\eta}\leq l(\xi\cdot \nu_E)\leq\|\xi_{|\Gamma}\cdot \nu_E\|_\infty\frac{l(\zeta_0)}{\eta},$$
which yields the result with
$k=l(\zeta_0)/\eta$.
\qed
\end{proof}

\begin{remark}\label{rem:reg}\rm The condition \eqref{H} requires some regularity. Indeed, we prove that if $\partial E$ has a cusp, \eqref{H} is not verified (see subsection (\ref{ex2}) below). However,

\begin{proposition}\label{Liprop}
If $E$ is a Lipschitz bounded open set, then $E$ verifies
 \eqref{H}.
 \end{proposition}

\begin{proof} Given $x\in\partial E$, up to a change of coordinates, we may assume that $\partial E$ is locally the graph of a Lipschitz function $f^x:\Omega_x\subset\R^{N-1}\rightarrow\R$
 on a neighborhood of $x\in \partial E$, and $E$ is locally below this graph. The normal $\nu_E$ to the boundary is given by $\nu_E(x',f^x(x'))=(-\nabla f^x(x'),1)/\sqrt{|\nabla f^x(x')|^2+1}$ so that, if $\xi_0(x',x_N)=(0,1)$, we have on $\partial E$:
 $$(\xi_0\cdot \nu_E)(x',x_N)=1/\sqrt{|\nabla f^x(x')|^2+1}\geq 1/\sqrt{\|\nabla f^x\|^2_\infty+1}=\eta>0.$$
Using the compactness of $\partial E$ and a regular partition of
unity, we construct a global vector field $\xi_0$ satisfying
(\ref{H}). \qed
\end{proof}\end{remark}

\subsection{A counterexample}\label{ex2} The following example shows that, when Lipschitz continuity of the boundary does not hold, not only condition (\ref{H}) may not be satisfied, but positivity of $l$ may not imply its $L^\infty$-continuity.

We consider again Example \ref{ex0} with a cusp.
We saw -see (\ref{P})- that
$$ \forall \xi\in \Theta,\; \mathcal{P}_E'(0)\cdot\xi=\int_\Gamma H_\Gamma \; \xi\cdot \nu_E \, d\mathcal{H}^1
-2\,\xi_2(0,0)=l(\xi_{|\Gamma}\cdot \nu_E).$$

Let us consider the shape functional
$$\mathcal{J}_E(\theta)=P(E_\theta)-\int_{E_\theta}\!\!H(z)dz,$$
where $H$ is the extension to $\R^2$ of $H_\Gamma$ given by
$$\forall (x,y)\in\R^2,\; H(x,y)=f''(x)[1+f'(x)^2]^{-3/2}=\alpha(\alpha-1)|x|^{\alpha-2}\big[1+\alpha^2|x|^{2\alpha-2}\big]^{-\frac{3}{2}}.$$
Then, $\frac{\partial H}{\partial x} \mathop{\sim}_{x\rightarrow 0} Cx^{-2\alpha}$ so that $H\in W^{1,1}_{loc}(\R^2)$ if $\alpha<1/2$ (which we assume).
This is enough to differentiate $\mathcal{J}_E$ and to obtain
$$\mathcal{J}_E'(0)\cdot\xi=-2\,\xi_2(0,0)=l_1(\xi_{|\Gamma}\cdot \nu_E),$$
where the existence of $l_1$ is given by Theorem \ref{the1}.

The distribution $\mathcal{J}_E'(0)$ is of order 0 (a Dirac mass)
and is \textbf{positive} since, according, for instance, to
(\ref{xi2})
$$ [\xi\cdot \nu_E\geq 0\;{\rm on}\;\partial E\setminus {(0,0)}]\;\Rightarrow\;-\xi_2(0,0)\geq 0\;.$$
However, $l_1$ {\em is not continuous for the uniform norm}. Indeed, let $\rho\in\mathcal{C}^\infty_0(\R^2), 0\leq\rho\leq 1$ with support in $B(0,\eta)$ and $\rho\equiv 1$ on $B(0,\eta/2)$. We choose $\xi=(0,1)\rho$. Then, $\xi_2(0,0)=1$, while $$\|\xi\cdot \nu_E\|_\infty\leq \sup_{|x|\leq \eta}[1+\alpha^2|x|^{2\alpha-2}]^{-1/2}=[1+\alpha^2\eta^{2\alpha-2}]^{-1/2}\to 0\;{\rm as}\;\eta\to 0,$$
so that there does not exist $k$ such that
$$\forall \xi\in\Theta, |\mathcal{J}_E'(0)\cdot\xi|=2|\xi_2(0,0)|\leq k\|\xi\cdot \nu_E\|_\infty.$$

\subsection{Perimeter of convex and pseudoconvex sets}
Let us first apply the results of the previous sections to convex sets.
\begin{proposition}\label{convexprop}
Let $C$ be a bounded convex set in $\R^N$. We let $\mathcal{P}_C(\theta):=P((Id+\theta)(C))$.
Then $\mathcal{P}_C'(0)$ is positive.
\end{proposition}
\begin{proof}
(See also \cite{BG} for a different approach). According to Lemma \ref{croiss}, it is sufficient to prove that, for small $\theta\in\Theta$,
$$[C\subset C_\theta]\;\Rightarrow\;[P(C)\leq P(C_\theta)]\;.$$
For this, we use the classical property (see the Appendix)
\begin{eqnarray}\label{capH}
\forall D\subset\R^N\;{\rm measurable,\;bounded}, \forall H\subset \R^N\;{\rm closed\;halfspace},\; P(D\cap H)\leq P(D) \;.
\end{eqnarray}
Next, we use that $\overline{C}=\cap_{n\geq 0}H_n$ where $H_n$ are closed halfspaces. For the sequence defined by: $\forall n\geq 0,D_{n+1}=D_n\cap H_n; D_0=C_\theta$, we have
$$P(D_{n+1})\leq P(D_n),\;P(C)=P(\overline{C})\leq\liminf_{n\to\infty}P(D_n)\leq P(C_\theta).$$\qed
\end{proof}
\noindent{\bf Proof of Theorem \ref{thm:convex}.} A bounded convex open set is Lipschitz. We apply Propositions \ref{convexprop},\ref{Liprop},\ref{prop:pos} to obtain the estimate (\ref{PC1}). Then (\ref{PC2}) follows as indicated in the remark following Theorem \ref{thm:convex}.\qed

The above results may be extended to pseudoconvex sets in the
spirit of \cite{BGM}: our approach allows us to assume that the
boundary is only Lipschitz (rather than
 $\mathcal{C}^{1,\alpha}$).

\begin{proposition}\label{pseudo} Let $E\subset\R^N$ be a Lipschitz bounded open set. Assume that $E$ is {\em pseudo-convex}, that is, with the notations of the proof of Proposition \ref{Liprop}: for all $x\in\partial E$ and for all $v:\Omega_x\to \R$ Lipschitz continuous, nonnegative and compactly
supported, one has
\begin{eqnarray}\label{pseudoing}\int_{\Omega_x}\sqrt{1+|\nabla f^x(x')|^2}\,dx'\leq \int_{\Omega_x}\sqrt{1+|\nabla [f^x+v](x')|^2}\,dx'.
\end{eqnarray}
Then the conclusions of Theorem \ref{thm:convex} hold.
\end{proposition}

\noindent {\bf Proof:} Like in the proof of Theorem \ref{thm:convex}, it is sufficient to prove that $\mathcal{P}_E'(0)\geq 0$. According to Lemma \ref{croiss} (and its proof) and by localization (using a partition of unity), it is sufficient to prove that, for $\eta$ and $t$ small
$$[\xi\cdot\nu_E\geq 0\;on \;\partial E, Supp(\xi)\subset B(x,\eta)]\Rightarrow [P(\gamma_t(E))\geq P(E)].$$
Like in the proof of Lemma \ref{croiss}, we have $E\subset
\gamma_t(E)$. But $\partial \gamma_t(E)$ is also the graph of a
Lipschitz function $v:\Omega_x\to [0,+\infty)$: indeed, if we
denote
$\gamma_t(x',f^x(x'))=\big(\gamma'(x',f^x(x')),\gamma_N(x',f^x(x'))\big)$,
then, $x'\to \gamma'(x',f^x(x'))$ is a local bi-Lipschitz
diffeomorphism since $D_x\gamma_t(x)$ is closed to the identity
($t$ small) and $f^x$ is Lipschitz continuous. Then the inequality
(\ref{pseudoing}) exactly says that $P(\gamma_t(E))\geq P(E)$.\qed

\subsection{Integrals on variable boundaries}
We now consider the more general shape functional
$\theta\to\mathcal{J}_E(\theta)=\int_{\Gamma_\theta}g\,
d\mathcal{H}^{N-1}$ where $\Gamma=\partial^*\!E$ and $g\in
\mathcal{C}^{1}(\R^N)$, with $(g,\nabla g)$ bounded. By change of
variable (see \cite{HP} for the regular case, and \cite{Giusti}
for the case of sets of finite perimeter we are using here), we
obtain :
$$\forall \theta \textrm{ such that }\|\theta\|_{1,\infty}<1,\;\;\mathcal{J}_E(\theta)=\int_{\Gamma}g\circ(Id+\theta)\;Jac_\Gamma(Id+\theta)
d\mathcal{H}^{N-1}$$ where $Jac_{\, \Gamma} T:= |det(DT)| \,.
\,|^tDT^{-1}\cdot\nu_E|$. We can prove (as in the regular case
done in \cite{HP}) that this functional is differentiable at
$\theta=0$. Using
$$Jac_{\, \Gamma} (Id+t\xi)=1+t\nabla_\Gamma\cdot \xi+o(t),\;\;g\circ(Id+t\xi)=g+t\nabla g\cdot\xi+o(t),$$
we obtain
$$\mathcal{J}_E'(0)\cdot\xi=\int_\Gamma \nabla g\cdot\xi
+g\nabla_\Gamma\cdot \xi\;.$$ We rewrite this with the help of the
formulas
$$\nabla g\cdot\xi+g\nabla_\Gamma\cdot \xi=\nabla g\cdot\xi+\nabla_\Gamma\cdot (g\xi)-\nabla_\Gamma g\cdot\xi=(\nabla g\cdot\nu_E)(\xi\cdot\nu_E)+\nabla_\Gamma\cdot (g\xi).$$
$$\mathcal{J}_E'(0)\cdot\xi=\int_\Gamma (\nabla g\cdot\nu_E)(\xi\cdot \nu_E)+\int_\Gamma \nabla_\Gamma\cdot (g\xi)
=\int_\Gamma (\nabla g\cdot\nu_E)(\xi\cdot
\nu_E)+\mathcal{P}'_E(0)\cdot(g\xi).$$ If $l_g, l_1$ denote the
linear functionals respectively associated to $\mathcal{J}_E'(0)$
and $\mathcal{P}'_E(0)$ by Theorem \ref{the1}, we obtain
$$\forall \varphi\in\Phi(\Theta),\;\;l_g(\varphi)=\int_\Gamma (\nabla g\cdot\nu_E)\,\varphi+l_1(g\varphi).$$
As a consequence, since $\nabla g$ is bounded, {\em $l_g$ is $L^\infty$-continuous if and only if $l_1$ is $L^\infty$-continuous.}
For example, it is the case when $E=C$ is a convex set, and then (see Theorem \ref{thm:convex}),
$$\mathcal{J}_C'(0)\cdot\xi=\int_{\partial C} (\nabla g\cdot\nu_C)(\xi\cdot \nu_C)+\int_{\partial C}g\,\xi\cdot d\overrightarrow{\textbf{H}}.$$

\subsection{Application to the Dirichlet energy}
Given $f\in L^2(\R^N)$, to each open bounded subset $E$ of $\R^N$,
we may associate the solution to the Dirichlet problem
\begin{eqnarray}\label{Dir}
u_E\in H^1_0(E),\;\;-\Delta u_E =f\;on\;E\;,
\end{eqnarray}
and its Dirichlet energy
\begin{eqnarray}\label{energy}
J(E)=\int_E\frac{1}{2}|\nabla u_E|^2-f\,u_E=-\frac{1}{2}\int_E|\nabla u_E|^2.
\end{eqnarray}
The last equality may be obtained by multiplying (\ref{Dir}) by
$u_E$ and integrating by parts. We also know that
\begin{eqnarray}\label{Dir1}
J(E)=\min\{\int_E\frac{1}{2}|\nabla v|^2-f\,v;\;v\in H^1_0(E)\}\;.
\end{eqnarray}

{\em For $E\subset\R^N$ given, let us analyze the derivative of}
$\theta\in \Theta \to \mathcal{J}_E(\theta):=J(E_\theta).$ \\

As proved for instance in \cite{HP}, for each bounded open subset
$E$, $[\theta\in\Theta\mapsto \mathcal{J}_E(\theta)]$ is
Fr\'echet-differentiable at $\theta=0$ and this does not require
any regularity of $E$.

 The structure of this derivative is given by Theorem \ref{the1} as soon as $\mathcal{J}_E(\theta)$ depends only on $\chi_{E_\theta}$, i.e., as soon as $J$ satisfies the property (\ref{hypJ}). This is the case if $E$ is regular since then $E_{\theta_1}=E_{\theta_2}$ a.e. implies $E_{\theta_1}=E_{\theta_2}$ everywhere (we use that $Id+\theta_1,Id+\theta_2$ are diffeomorphisms).

 Let us describe an optimal and (very) weak notion of regularity of $E$ for which (\ref{hypJ}) still holds for the Dirichlet energy. Let us first recall some definitions for $E\subset\R^N$ only measurable and bounded (see e.g. \cite{HP}):
  \begin{eqnarray}\label{h1}
 H^1_0(E):=\{v\in H^1(\R^N);\;v=0\;quasi-everywhere\; outside\;E\}\;.
 \end{eqnarray}
 \begin{eqnarray}\label{h2}
 \widehat{H}^1_0(E):=\{v\in H^1(\R^N);\;v=0\;a.e.\;outside\;E\}\;.
 \end{eqnarray}
It is well-known that the "right" extension for the definition of $H^1_0(E)$ in a non-regular setting is given by (\ref{h1}) where $"quasi-everywhere"$ means "everywhere except on a set of capacity zero" for the usual $H^1(\R^N)$-capacity. In particular, it coincides with the usual definition of $H^1_0(E)$ for any bounded open set $E$. The second definition (which is different and leads to a larger space) has also proven to be useful in shape optimization problems, due to the simplicity of its definition and due to the fact that we are then eventually lead to work with $E$'s which are regular enough to satisfy
\begin{eqnarray}\label{regular}
\widehat{H}^1_0(E)=H^1_0(E).
\end{eqnarray}
Here, we adopt this second point of view {\em and we assume that
$E$ satisfies (\ref{regular})}. Then for $\theta$ small in
$\Theta$, since $Id+\theta$ is a diffeomorphism, it is easy to
verify that $E_\theta$ also satisfies (\ref{regular}). It follows
that the Dirichlet energy $\theta\mapsto
\mathcal{J}_E(\theta)=J(E_\theta)$ where $J$ is defined according
to (\ref{energy}) with $E_\theta$ in place of $\Omega$ satisfies
the hypothesis (\ref{hypJ}). Indeed, we then have for small
$\theta_1,\theta_2\in\Theta$
$$[E_{\theta_1}=E_{\theta_2}\;a.e.]\Rightarrow[\widehat{H}^1_0(E_{\theta_1})=\widehat{H}^1_0(E_{\theta_2})]\Rightarrow [H^1_0(E_{\theta_1})= H^1_0(E_{\theta_2})]\Rightarrow[J(E_{\theta_1})=J(E_{\theta_2})].$$

And we have the following, in the spirit of Theorem \ref{the1}:
\begin{proposition} Let $E$ be measurable, bounded, with finite perimeter and satisfying (\ref{regular}). Then, there exists $l:\Phi(\Theta)\to\R$ linear continuous and {\bf positive} such that
$$\forall \xi\in\Theta,\;\;\mathcal{J}_E'(0)\cdot \xi=-l(\xi_{|\partial^*\!E}\cdot \nu_E).$$
Then, if moreover $E$ has a Lipschitz boundary, then $l$ is $L^\infty$-continuous.
\end{proposition}
\noindent{\bf Proof.} For the existence of $l$, we apply Theorem \ref{the1}. For the positivity of $l$, we apply Lemma \ref{croiss} after noticing that, thanks to (\ref{regular})
$$[E\subset E_\theta\;a.e.]\Rightarrow [H^1_0(E)\subset H^1_0(E_\theta)]\Rightarrow [J(E)\geq J(E_\theta)].$$
For the last remark, we apply Propositions \ref{Liprop}-\ref{prop:pos}.\qed

\noindent{\bf Remark:} When $E$ is regular enough, we have
$$-\mathcal{J}_E'(0)\cdot\xi=\frac{1}{2}\int_{\partial E}|\nabla u_E|^2(\xi\cdot \nu_E),$$
and the structure as well as the positivity of the derivative are then obvious. The above proposition says that positivity (in our sense) remains true without much regularity on $E$ and even if one is not allowed to write the above expression which, in particular, requires that a trace of $\nabla u_E$ be defined on the boundary. If $E$ has a Lipschitz boundary, the $L^\infty$-continuity of $l$ provides some kind of $L^1$-bound on the trace of $|\nabla u_E|^2$.\\

The same property holds for the derivative of {\em the eigenvalues of the Laplace operator with homogeneous Dirichlet boundary conditions}. They may be defined for measurable sets $\Omega$ satisfying (\ref{regular}). If we denote the $k$-th eigenvalue by $\lambda_k(\Omega)$, then we have also the monotonicity property: $[E\subset E_\theta]\Rightarrow[\lambda_k(E)\geq\lambda_k(E_\theta)].$
If they are simple eigenvalues, they are also Fr\'echet-differentiable on $\Theta$ at $\theta=0$ without any regularity on $E$. Their derivatives are given in the regular case by
$$-(\lambda_k)_E'(0)\cdot \xi=\int_{\partial E}[\nabla u_k\cdot \nu_E]^2(\xi\cdot \nu_E),$$
where $u_k$ is the associated eigenfunction. The positivity of the derivative is obvious on this formula. Using the monotonicity property of $\lambda_k$ and the same arguments as in the previous proposition, we prove that, even with poor regularity on $E$,
$(\lambda_k)_E'(0)\cdot \xi=-l(\xi_{|\partial^*\!E}\cdot \nu_E)$ for some $l\geq 0$.

\section{Appendix}

\noindent{\bf Some facts of geometric measure theory}

Let first recall some known facts from geometric measure theory (see e.g. \cite{EG} and \cite{Federer}). Given $E$ a measurable subset of $\R^N$, its perimeter is defined by
$$P(E):=\textrm{sup} \left\{
 \int_{E} \nabla\cdot  \varphi\,;\; \varphi=(\varphi_1,...,\varphi_N) \in {\mathcal C}_0^\infty(\R^N,\R^N),\;\sup_{x\in\R^N}\sum_{1\leq i\leq N}\varphi_i^2(x)\leq 1 \right\}.$$

According to the duality formula
$$\int_{E} \nabla\cdot \varphi\, dx =-\left<\nabla
\chi_E,\varphi\right>_{\mathcal{D}'(\R^N,\R^N)\times
\mathcal{C}_0^\infty(\R^N,\R^N)},$$ where
$\mathcal{D}'(\R^N,\R^N)$ denotes the space of $\R^N$-valued
distributions on $\R^N$, saying that $E$ is of finite perimeter
means that the distribution $\nabla \chi_E$ is a bounded
(vector-valued) measure on $\R^N$ and its total mass is given by
$P(E)$. This measure has the polar decomposition $\nabla
\chi_E=-|\nabla \chi_E|\nu^E$ where $|\nabla\chi_E|$ is a positive
measure, $\nu^E\in L^\infty(\R^N,|\nabla\chi_E|)$ and
$|\nu^E|=1\;\;|\nabla\chi_E|$-a.e..

Next, we may define a pointwise representation of $\nu^E$ together
with the reduced boundary of $E$ as follows.
\begin{defin}Let $E$ be a set of finite perimeter. Its reduced boundary $\partial^*\!E$ is defined as the subset of $x\in\R^N$ such that  $|\nabla\chi_E|(B(x,\rho))>0$ for all $\rho>0$ and
such that the limit
$$\nu_E(x):= -\lim_{\rho\rightarrow0}\frac{\nabla \chi_E(B(x,\rho))}{|\nabla \chi_E|(B(x,\rho))}$$
exists in $\R^N$ and its Euclidian norm equals $1$.
\end{defin}
Actually, $\nu_E$ is a pointwise representation of $\nu^E$ and may be interpreted as the generalized outward unit normal to $E$. More precisely, we have the following.
\begin{proposition}\label{nablachiE}
The measure $|\nabla\chi_E|$ is carried by $\partial^*\!E$ and is
the restriction to $\partial^*\!E$ of the $(N-1)$-Hausdorff measure
$\mathcal{H}^{N-1}$ on $\R^N$:
$$|\nabla \chi_E| = {\mathcal{H}^{N-1}}_{|\partial^*\!E},$$
$$\textrm{i.e., } \forall B \textrm{ Borel set},\; |\nabla \chi_E|(B) = \mathcal{H}^{N-1}(\partial^*\!E\cap B).$$
Moreover, $$\nabla\chi_E=-|\nabla\chi_E|\,\nu_E=-\nu_E\,{\mathcal{H}^{N-1}}_{|\partial^*\!E}.$$
\end{proposition}

\noindent{\bf A technical lemma}
\begin{lemma}\label{deriv}
Let $E\subset \R^N$ with finite perimeter, $\xi\in \Theta$ and
$\gamma_t$ its associated flow, that is, the solution of
$$\forall x\in\R^N,\;\forall t\in\R,\;\frac{d}{dt}\gamma_t(x)=\xi(\gamma_t(x)),\;\;\gamma_0(x)=x\,.$$
Then, in the sense of distributions in $\R\times\R^N$,
$$\frac{\partial}{\partial t}
(\chi_E\circ\gamma_t)=\nabla\chi_E\circ\gamma_t\cdot \frac{\partial}{\partial
t}\gamma_t =(\nabla\chi_E\cdot \xi)\circ\gamma_t,$$
which means that, for all $\varphi\in\mathcal{C}_0^\infty(\R\times\R^N)$,
$$-\int_{\R\times\R^N}\frac{\partial \varphi}{\partial t}\,\chi_E\circ\gamma_t=\int_{\R\times\R^N}\varphi\,d[(\nabla \chi_E\cdot\xi)\circ\gamma_t]dt=...$$
$$=\int_{\R\times\R^N}(\varphi\circ\gamma_t^{-1})\det D(\gamma_t^{-1})\xi\cdot d(\nabla\chi_E)dt=-\int_\R\int_{\partial^*\!E}(\varphi\circ\gamma_t^{-1})\det D(\gamma_t^{-1})(\xi\cdot\nu_E)\,d\mathcal{H}^{N-1}dt.$$
\end{lemma}

\noindent {\bf Proof:}
Let $f_n\in\mathcal{C}^{\infty}_0(\R^N)$ such that
$$f_n \stackrel{L^1}{\longrightarrow} \chi_E,\;\;
\nabla f_n \stackrel{(\mathcal{C}_b)'}{\longrightarrow}
\nabla \chi_E,$$
where $\mathcal{C}_b$ denotes the space of bounded continuous functions from $\R^N$ into $\R$.
Then $$\frac{\partial}{\partial
t}(f_n\circ\gamma_t)=\nabla
f_n\circ\gamma_t\cdot\frac{\partial}{\partial t}\gamma_t =(\nabla
f_n\circ\gamma_t)\cdot(\xi\circ\gamma_t)=(\nabla f_n\cdot\xi)\circ\gamma_t.$$
Using arbitrary $\varphi\in\mathcal{C}_0^\infty(\R\times\R^N)$, we may rewrite this as
$$-\int_{\R\times\R^N}\frac{\partial \varphi}{\partial t}\,f_n\circ\gamma_t=\int_{\R\times\R^N}\varphi\,(\nabla f_n\cdot\xi)\circ\gamma_t,$$
or, after change of variable
$$-\int_{\R\times\R^N}\frac{\partial \varphi}{\partial t}\circ\gamma_t^{-1}\,f_n\;\det D(\gamma_t^{-1})=\int_{\R\times\R^N}(\varphi\circ\gamma_t^{-1})(\nabla f_n\cdot\xi)\;\det D(\gamma_t^{-1}).$$
Since $\det D(\gamma_t^{-1})$ and $\xi$ are continuous and $\varphi$ is $\mathcal{C}^\infty$ and compactly supported, we may pass to the limit in this equality to obtain
$$-\int_{\R\times\R^N}\frac{\partial \varphi}{\partial t}\circ\gamma_t^{-1}\,\chi_E\;\det D(\gamma_t^{-1})=\int_{\R\times\R^N} (\varphi\circ\gamma^{-1}_t)\det D(\gamma_t^{-1})\,\xi\cdot d\nabla \chi_Edt.$$

The formula of Lemma \ref{deriv} follows, using also Proposition \ref{nablachiE}.\qed

\noindent {\bf A proof of (\ref{capH}):}\\
We choose coordinates $(x',x_N)\in\R^{N-1}\times\R$ so that $H=\{x_N\leq 0\}$. For $t\in (0,1]$, we set
$$f_t(x):=\left\{\begin{array}{ccl}x&if&x\in H\\
(x',tx_N)&if&x\in H^c\;.\end{array}\right.$$
Let us check that:
(i) $P(f_t(D))\leq P(D)$ and (ii) $\chi_{f_t(D)}$ tends a.e. to $\chi_{D\cap H}$ as $t$ tends to $0^+$.\\
Then $P(D\cap H)\leq P(D)$ will follow by the lower semi-continuity of $P(\cdot)$.
The point (ii) comes easily from the fact that $D$ is bounded. For (i), we first notice that,
$$|\nabla \chi_D|(H)=\int_{\partial^*\!D\cap H}d\mathcal{H}^{N-1}=\int_{\partial^*\!f_t(D)\cap H}d\mathcal{H}^{N-1}=|\nabla \chi_{f_t(D)}|(H).$$ Next, let $\varphi\in\mathcal{C}_0^\infty(H^c)$ with $\|\varphi(x)\|\leq 1$ for all $x\in H^c$. We have
$$\int_{f_t(D)\cap H^c}\nabla\cdot\varphi\,(x)\, dx=\int_{D\cap H^c}t\,\nabla\cdot\varphi\,(x',tx_N)\,dx=\int_{D\cap H^c}\nabla\cdot\psi\,(x',x_N)\,dx,$$
where we set
$\psi(x',x_N)=\big(t\varphi_1(x',tx_N),...,t\varphi_{N-1}(x',tx_N),\varphi_N(x',tx_N)\big)$.
As $\|\psi(x)\|\leq 1$, we deduce that
$$\int_{f_t(D)\cap H^c}\nabla\cdot\varphi\,(x)\, dx\leq |\nabla \chi_D|(H^c).$$
Since $\varphi$ is arbitrary in $\mathcal{C}_0^\infty(H^c)$ where $H^c$ is open, it follows that $|\nabla \chi_{f_t(D)}|(H^c)\leq|\nabla \chi_D|(H^c).$\qed

\bibliographystyle{plain}

\end{document}